\documentclass[12pt]{iopart}
\usepackage{amsfonts}
\usepackage[english]{babel}
\usepackage{euscript}
\newcommand{\textfrac}[2]{{\textstyle\frac{#1}{#2}}}
\newcommand{\kpp}{\kappa}
\newcommand{\w}{\wedge}
\newcommand{\qqq}{\quad\quad\quad}
\newcommand{\brk}{\]\[\fl\qqq}

\begin{document}

\title[Contact integrable extensions for r-mmdKP equation]{%
Contact Integrable Extensions of Symmetry Pseudo-Group and Coverings for the r-th Double
\\
Modified Dispersionless Kadomtsev--Petviashvili Equation
\footnote[7]{The work  was partially supported by the joint grant 09-01-92438-KE\_a of RFBR (Russia) and Consortium E.I.N.S.T.E.IN (Italy).}
}

\author{Oleg I. Morozov}

\address{Department of Mathematics, Moscow State Technical University
\\
of Civil Aviation, Kronshtadtskiy Blvd 20, Moscow 125993, Russia
\\
oim{\symbol{64}}foxcub.org}

\begin{abstract}
We find contact integrable extensions and coverings for the r-th double
modified dis\-per\-si\-on\-less Kadomtsev--Petviashvili equation.
\end{abstract}


\vskip 40 pt

We consider the r-th double modified dispersionless Kadomtsev--Petviashvili equa\-ti\-on
\begin{eqnarray}
u_{yy} &=& u_{tx}
+ \left(\frac{(\kappa+1)\,u_y^2}{u_x^2}-\frac{u_t}{u_x}+\kappa\,u_x^\kappa u_y
+\frac{(\kappa+1)^2}{2\,\kappa+3}\,u_x^{2(\kappa+1)}
\right)\,u_{xx}
\nonumber
\\
&&
-\kappa\,\left(\frac{u_y}{u_x}+u_x^{\kappa+1}\right)\,u_{xy}
\label{rmmdKP}
\end{eqnarray}
with $\kappa \not \in \{-2, -3/2, -1\}$. This equation appears from the differential covering,
\cite{KV84,KLV,KV89},
\begin{equation}
\left\{
\begin{array}{l}
u_t =
\displaystyle{ 
\left(\frac{(\kappa+2)^2}{2\,\kappa+3}\,u_x^{2(\kappa+1)}
+(\kappa+2)\,w_x\,u_x^{\kappa+1}+\frac{\kappa+1}{2} \,w_x^2-w_y\right) u_x}
\\
u_y = -\left(u_x^{\kappa+1}+w_x\right) u_x
\end{array}
\right.
\label{covering_over_rmdKP}
\end{equation}
over the r-th modified dispersionless Kadomtsev--Petviashvili equation, \cite{Blaszak},
\begin{equation}
w_{yy} = w_{tx}+\left(\textfrac{1}{2}\,(\kappa+1)\,w_x^2+w_y\right)\,w_{xx} +\kappa\,w_x\,w_{xy},
\label{rmdKP}
\end{equation}
see
\cite{ChangTu},
\cite{KonopelchenkoAlonso},
\cite{Pavlov2006},
\cite{Morozov2010}.
Namely, excluding $w$ from (\ref{covering_over_rmdKP})  yields Eq. (\ref{rmmdKP}).

We apply the method of contact integrable extensions, \cite{Morozov2009}, to find differential co\-ve\-rings of Eq. (\ref{rmmdKP}). The method starts from computing Maurer--Cartan forms and struc\-tu\-re equations for the symmetry pseudo-group via approach of 
\cite{Morozov2002a,Morozov2004b}. The structure equations read
\[
\fl
d\theta_0=
(
\theta_{22}
+
(U_2
-(\kpp+1)^2\,(U_1-2\,(\kpp+1))\,\xi^1
-(U_4
+(\kpp+1)(\kpp+2)\, U_3
-\kpp\,(\kpp+1)\,U_1
\brk    
+(\kpp+1)^2(\kpp+2)(\kpp^2+6\,\kpp+4))\,(\kpp+1)^{-2}(\kpp+2)^2\,\xi^2
-(
U_4
-\kpp\,(\kpp+1)\,U_1
\brk    
+(\kpp+1)^2(\kpp+2)\,(2\kpp+3))\,(\kpp+1)^{-1}\,(\kpp+2)^{-1}\,\xi^3)
\w \theta_0
+\xi^1\w \theta_1 +\xi^2\w \theta_2
\brk    
+\xi^3\w \theta_3, 
\]
\[
\fl
d\theta_1
=
(\kpp+1)\,(2\,\theta_1+(\kpp+1)(\kpp+2)^2\,\theta_2
-(\kpp+2)\,\theta_3))\w \theta_0 +\xi^1\w \theta_{11}+\xi^3 \w \theta_{13}
\brk    
+(\kpp+1)(\kpp+2)\,\theta_2 \w \theta_3
-(2\,(\kpp+1)\,(\theta_2-2\,(U_1-2\,(\kpp+1)(\kpp+2))\,\xi^1+\xi^2)
\brk    
+((2\,\kpp+3)\,U_1-(\kpp+1)(\kpp+2)(3\,\kpp + 4))(\kpp+2)^{-1}\,\xi^3)\w \theta_1
\brk    
+\xi^2\w (U_3\,\theta_0+U_1\,\theta_3+\theta_{12}),
\]
\[
\fl
d\theta_2=
\theta_0\w\theta_{22}
+\left((\kpp+1)\,(U_1-2\,(\kpp+1)(\kpp+2))\,\xi^1+\textfrac{1}{2}\,U_1\,\xi^3\right)\w \theta_2
+\xi^1\w \theta_{12}+\xi^2\w \theta_{22}
\brk    
+\xi^3\w \theta_{23},
\]
\[
\fl
d\theta_3=
(
(\kpp+1)\,(\theta_3
-(\kpp+1)(\kpp+2)\,\theta_2
+(\kpp+1)(\kpp+2)(
(\kpp+3)\,U_1-(\kpp+2)\,(U_2+2))\,\xi^1)
\brk    
+U_3\,\xi^2+(U_2+U_4)\,\xi^3-(\kpp+1)(\kpp+2)\,\theta_{22})\w \theta_0
+\left((\kpp+1)\,\theta_3+\textfrac{1}{2}\,U_1\,\xi^2\right)\w \theta_2
\brk    
+(\kpp+1)(
2\,(U_1-2\,(\kpp+1)(\kpp+2))\,(\xi^1+(\kpp+2)^{-1}\,\xi^3)
-\xi^2)\w\theta_3
+\xi^1\w \theta_{13}
\brk    
+\xi^2\w \theta_{23}+\xi^3\w \theta_{12},
\]
\[
\fl
d\theta_{11}
=
\eta_1\w\xi^2+\eta_2\xi^3+\eta_3\w\xi^1
+
(
(4\,U_4-(\kpp+1)(\kpp-2)\,U_1-\kpp(\kpp+1)^2(\kpp^2-4)
\brk    
+(2\kpp+1)\,U_2)\,\theta_1
-(\kpp+1)^2(\kpp+2)\,(U_1-2\,(\kpp+1)(\kpp+2))\,(\theta_2+(\kpp+1)(\kpp+2)\,\theta_3)
\brk    
+(\kpp^2-1)(\kpp+2)\,\theta_{13}
+(\kpp+1)(\kpp\,(2\,U_5+3\,(\kpp+2)\,U_2)-U_1U_2
\brk    
-(\kpp+1)^2(\kpp+2)(3\,\kpp-2)((\kpp+3)\,\,U_1-2\,(\kpp+1)(\kpp+2))
)\,\xi^2
+((\kpp+1)(\kpp\,U_1
\brk    
-(\kpp+2)\,U_3+(\kpp+1)(\kpp+2)(3\kpp^2+6\kpp+4))-U_4)(\kpp+1)^{-2}(\kpp+2)^{-2}\,\theta_{11}
)\w\theta_0
\brk    
+(
(\kpp+1)(4\,U_1+(\kpp+1)(\kpp+2)\,(11\kpp+14))\,\theta_2
-(\kpp+1)(2\,U_1-(\kpp+1)(\kpp+2))\,\theta_3
\brk    
-(2\kpp+1)\,\theta_{12}
+4(\kpp+1)(\kpp+2)\,\theta_{23}
-(4\,U_4
-(\kpp+1)(2\,\kpp^2+3\,\kpp-4)\,U_1
\brk    
+(2\kpp+1)\,U_2
+4\,\kpp(\kpp+1)^2(\kpp+2)
-(2\kpp+3)(\kpp+2)^{-1}\,U_1^2
)\,\xi^2
\brk    
-(2(2\kpp+3)(\kpp+2)^{-1}\,U_5+(3\kpp+2)(\kpp+1)\,(U_2-(\kpp+3)(\kpp+1)\,U_1
\brk    
+2(\kpp+2)(\kpp+1)^2)\,\xi^3
)\w\theta_1
+(
(\kpp+2)(\kpp+1)^2\,(U_1-2\,(\kpp+1)\,(\kpp+2))\,\theta_3
\brk    
+(4\,\kpp+5)\,\theta_{11}
-(\kpp+1)(\kpp+2)\,\theta_{13}
)\w\theta_2
+((\kpp+2)\,\theta_{13}
-(
2\,U_5
\brk    
+(\kpp+1)(U_1^2-2(\kpp+2)(U_2+(\kpp+1)(\kpp+4)\,U_1+2(\kpp+1)^2(\kpp+2)))
)\,\xi^2
)\w\theta_3
\brk    
-(\theta_{22}
-((\kpp+1)(\kpp+9)\,U_1-U_2-14(\kpp+2)\,(\kpp+1)^2)\,\xi^1
-((\kpp+1)(\kpp\,U_1
\brk    
-(\kpp+2)\,(U_3-(3\kpp^2+6\kpp+4)(\kpp+1)))-U_4)(\kpp+1)^{-2}(\kpp+2)^{-2}\,\xi^2
\brk    
+(3\,(U_1-(\kpp+1)(\kpp+2))+(\kpp+1)^{-1}(\kpp+2)^{-1}\,U_4)\,\xi^3
)\w\theta_{11}
\brk    
+(
(\kpp+1)(U_1-2\,(\kpp+1)(\kpp+2))\,\theta_{12}-2\,U_1\,\theta_{13})\w\xi^2
\]
\[
\fl
d\theta_{12}= 
\eta_1\w\xi^1+\eta_4\w(\theta_0+\xi^2)
+\eta_7\w\xi^3
+(
\,(U_1-(\kpp+4)(\kpp+2)(\kpp+1))\,\theta_{22}
\brk    
+(\kpp+1)((\kpp+2)\,\theta_{23}
)
+(U_4+\textfrac{1}{2}\,(\kpp+1)(\kpp^2+2\kpp+2)\,U_1
+\kpp(\kpp+1)^2(\kpp+2))\,\theta_2
\brk    
+((\kpp+1)\,(\kpp\,U_1-(\kpp+2)\,(U_3-\kpp^2(\kpp+1)))-U_4)(\kpp+1)^{-2}(\kpp+2)^{-2}\,\theta_{12}
)\w\theta_1
\brk    
+((2\,\kpp+3)\,\theta_{12}
-(\kpp+1)(\kpp+2)\,\theta_{23}
-\textfrac{1}{2}\,U_1\,(\kpp+2)\,\theta_3
-U_5\,\xi^3
+(\textfrac{1}{2}\,U_1^2-(\kpp+1)\,U_1
\brk    
-U_4
-\kpp(\kpp+2)(\kpp+1)^2
)\,\xi^2
)\w\theta_2
+(\kpp+2)\,(\theta_{23}-(\kpp+1)\,\theta_{22})\w\theta_3
\brk    
-(\theta_{22} 
-((\kpp+1)(\kpp+5)\,U_1-U_2 -6(\kpp+2)(\kpp+1)^2 )\,\xi^1
+((\kpp+1)(\kpp\,U_1
\brk    
-(\kpp+2)(U_3-\kpp^2(\kpp+1))-U_4)(\kpp+2)^{-2}(\kpp+1)^{-2}\,\xi^2
-\textfrac{1}{2}\,
((3\kpp+8)(\kpp+1)\,U_1
\brk    
+2\,U_4
-4(\kpp+2)(\kpp+1)^2)(\kpp+1)^{-1}(\kpp+2)^{-1}\,\xi^3
)\w\theta_{12}
\brk    
+2(\kpp+1)(U_1-2\,(\kpp+2)(\kpp+1))\,(\theta_{22}\w\xi^2+\theta_{23}\w\xi^3)
-U_1\,\theta_{23}\w\xi^2
\]
\[
\fl
d\theta_{13}=
\eta_1\w\xi^3+\eta_2\w\xi^1+\eta_7\w\xi^2
+((\kpp+1)^2 (\kpp+2)^2\,\theta_{23}-U_3\,\theta_1
-(\kpp+1)(\kpp+2)\,\theta_{12}
\brk    
-\textfrac{1}{2}\, (\kpp+1)(\kpp+2)((\kpp+4)  (\kpp+1)\, U_1-4\,(U_2+ (\kpp+2) (\kpp+1)))\,\theta_2
\brk    
-((\kpp^2-1)\,U_1+U_2-U_4-2 (\kpp+2) (2 \kpp+1) (\kpp+1)^2)\,\theta_3
+((\kpp+1)\,(\kpp\,  U_1
\brk    
-(\kpp+2)\,(U_3-(2 \kpp^2+3 \kpp+2) (\kpp+1))-U_4)(\kpp+1)^{-2}(\kpp+2)^{-2}\,\theta_{13}
)\w\theta_0
\brk    
+(
\theta_{23}-(\kpp+1)(\kpp+2)\,\theta_{22}
+(U_2+U_4)\,\xi^3
+\textfrac{1}{2}\,(U_1-2(\kpp+2) (\kpp+1)^2)\,\theta_2
\brk    
+U_3\,\xi^2
)\w\theta_1
+(
(3\kpp+4)\,\theta_{13}
-(\kpp+1)(\kpp+2)\,\theta_{12}
-U_5\,\xi^2
-\textfrac{1}{2}\,(\kpp+1)((3 \kpp+4)\, U_1
\brk    
-4\, (\kpp+1) (\kpp+2))\,\theta_3
)\w\theta_2
+(
\theta_{12}+(\kpp+1)(\kpp+2)\,\theta_{23}
+((\kpp+1)(2\, U_1^2
\brk    
+(\kpp+2)((\kpp^2-\kpp-4)\,U_1-U_2-2 \kpp (\kpp+1)(\kpp+2)))-2 (\kpp+2) U_4)(\kpp+2)^{-1}\,\xi^2
\brk    
+((\kpp+1) (\kpp+2)(5 \kpp+2) ((\kpp+1)(\kpp+3)\,U_1- U_2-2 (\kpp+2) (\kpp+1)^2)
\brk    
-4 (\kpp+1) U_5)(\kpp+2)^{-1}\,\xi^3
)\w\theta_3
+\textfrac{3}{2}\,U_1\,\xi^2\w \theta_{12} 
-(\theta_{22}
-((\kpp+1)((\kpp+7)\,  U_1
\brk    
-10 (\kpp+2) (\kpp+1))-U_2)\,\xi^1
+((\kpp+1)(\kpp  U_1-(\kpp+2)  (U_3-(2 \kpp^2+3 \kpp+2) (\kpp+1))
\brk    
-U_4)(\kpp+1)^{-2}(\kpp+2)^{-2}\,\xi^2
-((3 \kpp+5) (\kpp+1) U_1-2 (\kpp+2) (2 \kpp+3) (\kpp+1)^2
\brk    
+U_4)(\kpp+1)^{-1}(\kpp+2)^{-1}\,\xi^3
)\w\theta_{13}
+(\kpp+1)(U_1-2 (\kpp+1) (\kpp+2))\,\theta_{23}\w\xi^2
\]
\[
\fl
d\theta_{22}=
\eta_4\w\xi^1+\eta_5 \w (\theta_0+\xi^2)
+\eta_6\w\xi^3
-(U_4-(\kpp+1)(\kpp  \,U_1+(\kpp+2)\, U_3
\brk    
+(\kpp+1)(\kpp+2) (\kpp^2+6 \kpp+4))
)(\kpp+1)^{-2}(\kpp+2)^{-2}\,
\theta_{22}\w (\theta_0 +\xi^2)
\brk    
+(((\kpp+1)^2 (U_1+2\kpp+4))-U_2)\,\xi^1
-((\kpp+1)(\kpp \, U_1+(\kpp+1)(\kpp+2) (3 \kpp+2))
\brk    
-U_4)(\kpp+1)^{-1}(\kpp+2)^{-1}\,\xi^3 
)\w\theta_{22}
\]
\[
\fl
d\theta_{23}= 
\eta_4\w\xi^3+\eta_6\w(\theta_0+\xi^2)+\eta_7\w\xi^1
+\textfrac{1}{2}\,(
U_1\,\theta_{22}
+(\kpp\,U_4-(\kpp+1)(\kpp^2 \, U_1
\brk    
-2  (\kpp+2)\, U_3+\kpp(\kpp+2) (3\kpp+2)(\kpp+1)^2)(\kpp+1)^{-1}(\kpp+2)^{-2}\,\theta_2
-2\,(\kpp  (\kpp+1)  U_1^2
\brk    
+((\kpp+2)\,U_3+U_4+(\kpp+2)(\kpp^2-3  \kpp-2)  (\kpp+1)^2)\,U_1
+\kpp (\kpp+1)  (\kpp+2)(U_4
\brk    
+(\kpp+2)\,  U_3))(\kpp+2)^{-2}\,\xi^3
-2\,(U_4-(\kpp+1)(\kpp\,  U_1+(\kpp+2) \, U_3
\brk    
+(\kpp+1)(\kpp+2)  (3\kpp+2)))(\kpp+1)^{-2}(\kpp+2)^{-2}\,\theta_{23})\w\theta_0 
+
\textfrac{1}{4}\,(
6\,(\kpp+2)\,\theta_{23}
\brk    
-8\,(\kpp+1)(\kpp+2)\,\theta_{22}
-2\,(\kpp\,U_4-(\kpp+1)(\kpp^2 \, U_1-2  (\kpp+2)\,U_3
\brk    
+\kpp(\kpp+1)(\kpp+2)(3\kpp+2))((\kpp+1)^{-1}(\kpp+2)^{-1}\,\xi^2
+(\kpp\,U_1^2-(\kpp+2)\,(2(\kpp+1) (\kpp^2
\brk    
+4  \kpp+2)\,U_1+2  (\kpp+4)\,U_2+4\,U_4))(\kpp+2)^{-1}\,\xi^3
)\w\theta_2
+(\kpp  (\kpp+1)  U_1+(\kpp+2)  U_3
\brk    
-U_4-(\kpp+2)  (3  \kpp+2)  (\kpp+1)^2)(\kpp+2)^{-1}\,\theta_3\w\xi^3
+(
\theta_{23}
-\textfrac{1}{2}\,U_1\,\xi^2
-2  (\kpp+1)  (U_1
\brk    
-2  (\kpp+1)  (\kpp+2))\,\xi^3
)\w\theta_{22}
+(((\kpp+3)  (\kpp+1)  U_1-U_2-2  (\kpp+2)  (\kpp+1)^2)\,\xi^1
\brk    
+(U_4-(\kpp+1) (\kpp\, U_1+(\kpp+2)\,U_3+(\kpp+2)  (3  \kpp+2)  (\kpp+1))(\kpp+1)^{-2}(\kpp+2)^{-2}\,\xi^2
\brk    
+\textfrac{1}{2}\,
(3  (\kpp+1)  (\kpp+2)  U_1+2  U_4-2  (\kpp+1)^2  (\kpp+2)^2
)(\kpp+1)^{-1}(\kpp+2)^{-1}\,\xi^3 )
\w\theta_{23}
\]
\[
\fl
d\xi^1
=
(
\theta_{22} 
+(2\,\kpp+3)\,\theta_2 
-((\kpp+1)((\kpp+3)\,U_1-2\,(\kpp+2))+U_4)\,((\kpp+1)(\kpp+2))^{-1}\xi^3
\brk    
+((\kpp+1)\,(\kpp\,U_1-(\kpp+2)\,(U_3-\kpp^2(\kpp+1)))-U_4)\,((\kpp+1)(\kpp+2))^{-2}\, (\theta_0+\xi^2)
) 
\w  \xi^1,
\]
\[
\fl
d\xi^2
=
(\theta_1-(\kpp+1)(\kpp+2)\,\theta_3-(\kpp+1)^2(\kpp+2)^2\,\theta_0) \w\xi^1
+(\theta_2+\theta_{22}
+(U_4-\kpp(\kpp+1)\,U_1
\brk    
-(\kpp+1)^2(\kpp+2)(3\kpp+2))(\kpp+1)^{-1}(\kpp+2)^{-1}\,\xi^3
-(U_4+(\kpp+1)(\kpp+2)\,U_3
\brk    
-\kpp(\kpp+1)\,U_1(\kpp+1)^2(\kpp+2)(\kpp^2+6\,\kpp+4))(\kpp+1)^{-2}(\kpp+2)^{-2}\,\theta_0
\brk    
+(U_2-(\kpp+1)(\kpp+2)\,U_1)\,\xi^1
)\w\xi^2
+(\theta_3 - (\kpp+1)(\kpp+2)\,\theta_3)\w\xi^3
\]
\[
\fl
d\xi^3
=
(
(\kpp+2)\,((\kpp+1)(\kpp\,\theta_0-2\,\theta_2)+\theta_3)
+(\kpp(\kpp+4)\,U_1-U_2-2\,(\kpp+1)^2(\kpp+2))\,\xi^3
\brk    
-U_1\,\xi^2
)\w\xi^1
+(
\theta_{22}+(\kpp+2)\,\theta_2
-(U_4+(\kpp+1)(\kpp+2)\,U_3-\kpp(\kpp+1)\,U_1
\]
\begin{equation}
\fl
\qqq
+(\kpp+1)^2(\kpp+2)(3\kpp+2))(\kpp+1)^{-2}(\kpp+2)^{-2})\,(\theta_0+\xi^2)
)\w\xi^3
\label{SEs_of_rmmdKP}
\end{equation}
The Maurer-Cartan forms $\theta_0$, ... , $\theta_{23}$, $\xi^1$, $\xi^2$, $\xi^3$ are 
\[
\fl
\theta_0 =u_{xx}u_x^{-2}\,(du-u_t\,dt-u_x\,dx-u_y\,dy)
\]
\[
\fl
\theta_1 = u_x^{-2\kappa-3}\,(du_t - u_{tt}\,dt-u_{tx}\,dx-u_{ty}\,dy)
-(\kappa+2)\,(u_y\,u_x^{-\kappa-2}-1)\,\theta_3
\brk    
+((\kpp+1)(\kpp+2)\,(u_yu_x^{-\kpp-2}-(2\kpp+3)^{-1})-u_tu_x^{-2\kpp-3} )
)\,\theta_2
\brk    
+(
u_tu_x^{-2\kpp-3}
+(\kpp+1)^2(\kpp+2)\,(u_yu_x^{-\kpp-2}-(2\kpp+5)(2\kpp+3)^{-1}
)\,\theta_0
\]
\[
\fl
\theta_2 = u_x^{-1}\,(du_x-u_{tx}\,dt-u_{xx}\,dx-u_{xy}\,dy)
\]
\[
\fl
\theta_3=u_x^{-\kappa-2}\,\left(du_y  - u_{tx}\,dt-u_{xy}\,dx
-E\,dy\right)
-(u_y\,u_x^{-\kappa-3}-\kappa-1)\,\theta_2
\brk    
-(u_y\,u_x^{-\kappa-3}+(\kappa+1)^2)\,\theta_0
\]
\[
\fl
\theta_{11}=u_{xx}^{-1}u_x^{-4\kpp-4}\,(du_{tt}-u_{ttt}\,dt-u_{ttx}\,dx-u_{tty}\,dy)
-2\,(\kpp+2)(u_yu_x^{-\kpp-2}-1)\,\theta_{13}
\brk    
-(2\,u_tu_x^{-2\kpp-3}-(\kpp+2)((\kpp+2)\,u_y^2u_x^{-2\kpp-4}-(2\kpp+3)\,u_yu_x^{-\kpp-2}
\brk    
+(2\kpp^2+9\kpp+8)(2\kpp+3)^{-1}))\,\theta_{12}
+A_{110}\,\theta_0+A_{111}\,\theta_1+A_{112}\,\theta_2+A_{113}\,\theta_3
\brk    
-(u_t^2u_x^{-4\kpp-6}
+(\kpp+1)^2(\kpp+2)^2\,(u_yu_x^{-\kpp-2}-(2\kpp+3)^{-1})^2
\brk    
+2\,(\kpp+1)(\kpp+2)(2\kpp+3)^{-1}\,u_x^{-2\kpp-3}
)\,\theta_{22}
-2\,(\kpp+2)\,((u_yu_x^{-\kpp-2}-1)u_tu_x^{-2\kpp-3}
\brk    
-(\kpp+1)(\kpp+2)(u_y^2u_x^{-2\kpp-4}+2\,(\kpp+2)(2\kpp+3)^{-3}\,u_x^{-\kpp-2}-2\kpp-3))\,\theta_{23}
\]
\[
\fl
\theta_{12}=u_{xx}^{-1}u_x^{-2\kpp-2}\,(du_{tx}-u_{ttx}\,dt-u_{txx}\,dx-u_{txy}\,dy)
-(u_yu_x^{-\kpp-2}-1)\,\theta_{23}
\brk    
+((\kpp+1)(\kpp+2)(u_yu_x^{-\kpp-2}-(2\kpp+3)^{-1})-u_tu_x^{-2\kpp-3})\,(\theta_{22}+\theta_3)
-\theta_1
\brk    
-(u_{txx}u_{xx}^{-2}u_x^{-2\kpp-1}+2\,u_{tx}u_{xx}^{-1}u_x^{-2\kpp-2}
-\textfrac{1}{2}\,(\kpp+2)\,(u_{xy}u_{xx}^{-1}\,((\kpp+2)u_yu_x^{-2\kpp-3}+\kpp\,u_x^{-\kpp-1})
\brk    
-(\kpp+2)\,(u_y^2u_x^{-2\kpp-4}-(\kpp+1)\,u_yu_x^{-\kpp-2})
+\kpp\,(\kpp+1))
)\,\theta_0
\]
\[
\fl
\theta_{13}= u_{xx}^{-1}u_x^{-3\kpp-3}\,(du_{ty}-u_{tty}\,dt-u_{txy}\,dx-\bar{\mathbb{D}}_t(E)\,dy)
-(2\kpp+3)(u_yu_x^{-\kpp-2}-1)\,\theta_{12}
\brk    
-(u_yu_x^{-\kpp-2}+(\kpp+1)^2)\,\theta_1-((u_yu_x^{-3\kpp-5}-(\kpp+1)\,u_x^{-2\kpp-3})\,u_t-(\kpp+1)\,(u_y^2u_x^{-2\kpp-4}
\brk    
-(\kpp+2)(2\kpp+3)^{-1}((2\kpp^2+5\kpp+4)u_yu_x^{-\kpp-2}-\kpp-1))\,\theta_{22}
+A_{130}\,\theta_0+A_{132}\,\theta_2
\brk    
+A_{133}\,\theta_3-(u_tu_x^{-2\kpp-3}-(\kpp+2)\,(u_y^2u_x^{-2\kpp-4}-(2\kpp+3)\,u_yu_x^{-\kpp-2}
\brk    
+2\,(\kpp+1)(\kpp+2)(2\kpp+3)^{-1}))\,\theta_{23}
\]
\[
\fl
\theta_{22}=u_{xx}^{-1}\,(du_{xx}-u_{txx}\,dt-u_{xxx}\,dx-u_{xxy}\,dy)-2\,\theta_2-u_xu_{xxx}u_{xx}^{-2}\,\theta_0
\]
\[
\fl
\theta_{23}=u_x^{-\kpp-1}u_{xx}^{-1}\,(du_{xy}-u_{txy}\,dt-u_{xxy}\,dx-\bar{\mathbb{D}}_x(E)\,dy)
-(u_yu_x^{-\kpp-2}\-\kpp-1)\,\theta_{22}-\theta_3
\brk    
+\textfrac{1}{2}\,(\kpp\,u_{xy}u_{xx}^{-1}u_x^{-\kappa-1}
-(\kappa+4)\,u_yu_x^{-\kappa-2}-\kappa\,(\kappa+1))\,\theta_2
\brk    
-(u_{xxy}u_{xx}^{-2}\,u_x^{-\kappa}-u_{xy}u_{xx}^{-1}\,u_x^{-\kappa-1}+u_yu_x^{-\kappa-2}+(\kappa+1)^2)\,\theta_0
\]
\[
\fl
\xi^1 = u_{xx}u_x^{2\kpp+1}\,dt
\]
\[
\fl
\xi^2 = u_{xx}u_x^{-1}\,dx
+(u_tu_x^{-\kpp-3}+(\kpp+2)\,(u_y^2u_x^{-2\kpp-4}-u_yu_x^{-2}+2\,(\kpp+1)^2(2\kpp+3)^{-1})
)\,\xi^1
\brk    
+(u_yu_x^{\kpp-2}-\kpp-1)\,\xi^3
\]
\begin{equation}
\fl
\xi^3 = u_{xx}u_x^{\kpp}\,dy+(\kpp+2)\,(u_yu_x^{-\kpp-2}-1)\,\xi^1,
\label{MCFs_of_rmmdKP}
\end{equation}
where $E$ is the right-hand side of Eq. (\ref{rmmdKP}), $\bar{\mathbb{D}}_t$, $\bar{\mathbb{D}}_x$ are restrictions of the total de\-ri\-va\-ti\-ves on Eq. (\ref{rmmdKP}), and  $A_{110}$, $A_{111}$, $A_{112}$, $A_{113}$, $A_{130}$, $A_{132}$, $A_{133}$ are functions of de\-ri\-va\-ti\-ves  of  $u$ of the first and the second orders. These functions are too long to write them in full. 
The forms $\eta_1$, ... , $\eta_7$ can be expressed from Eqs. (\ref{MCFs_of_rmmdKP}), (\ref{SEs_of_rmmdKP}).
The coefficients of the struc\-tu\-re equations depend on the invariants
\[
\fl
U_1 = (\kpp+2)\,(u_yu_x^{-\kpp-2}-u_{xy}u_{xx}^{-1}u_x^{-\kpp-1}+\kpp+1)
\]
\[
\fl
U_2  =u_{txx}u_{xx}^{-2}u_x^{-2\kpp-1}-(\kpp+2)\,u_{xxy}u_{xx}^{-2}u_x^{-\kpp}(u_yu_x^{-\kpp-2}-1)
-2\,u_{tx}u_{xx}^{-2}u_x^{-2\kpp-2}
+2\,u_tu_x^{-2\kpp-3}
\brk    
-(2u_yu_x^{-\kpp-2}-(\kpp+1)(\kpp+2))\,U_1
+2\,(\kpp+1)(\kpp+2)\,u_yu_x^{-\kpp-2}
\brk    
-u_{xxx}u_{xx}^{-2}\,(u_tu_x^{-2\kpp-2}-(\kpp+2)u_yu_x^{-\kpp-1}(u_yu_x^{-\kpp-2}-1)
\brk    
-2\,(\kpp+1)^2(\kpp+2)(2\kpp+3)^{-1}u_x)
+2\,(\kpp+1)(\kpp+2)(2\kpp^2+\kpp-2)(2\kpp+3)^{-1}
\]
\[
\fl
U_3 =u_{xxy}u_{xx}^{-2}u_x^{-\kpp}-u_{xxx}u_{xx}^{-2}u_x\,(u_yu_x^{-\kpp-2}+(\kpp+1)^2)+2\,(\kpp+2)^{-1}\,U_1
\brk    
-(\kpp+1)(\kpp^2+\kpp+2)
\]
\[
\fl
U_4 = (\kpp+1)\,(\kpp\,U_1-(\kpp+2)\,(U_3-(\kpp+1)\,(u_{xxx}u_{xx}^{-2}u_x+\kpp^2+5\kpp+2)))
\]
\[
\fl
U_5 = \textfrac{1}{2}\,(
(\kpp+2)u_{xx}^{-2}u_x^{-3\kpp-3}(u_xu_{txy}-u_{ty})
+u_{tx}u_{xx}^{-2}u_x^{-2\kpp-2}\,((\kpp+3)\,U_1-(\kpp+2)\,(u_yu_x^{-\kpp-2}
\brk    
+(\kpp+1)(\kpp+3)))
+((2\kpp+3)\,u_yu_x^{-\kpp-2}-1)\,U_1^2
-(\kpp+2)((\kpp+3)\,u_yu_x^{-\kpp-2}
\brk    
+2\kpp+1)\,U_2
-(
(\kpp+1)^{-1}u_tu_x^{-2\kpp-3}(\kpp(\kpp+1)^{-1}u_yu_x^{-\kpp-2}+2\kpp^2+5\kpp+4)
\brk    
+(\kpp+1)\,(\kpp\,u_y^2u_x^{-2\kpp-4}
-(2\kpp+3)^{-1}((2\kpp^4+9\kpp^3+7\kpp^2-13\kpp-18)\,u_yu_x^{-\kpp-2}
\brk    
-2\,(2\kpp^4+45\kpp^3+42\kpp^2+53\kpp+25))
)
)\,U_1
+(
u_tu_x^{-2\kpp-3}((\kpp+1)^{-1}u_yu_x^{-\kpp-2}-1)
\brk    
-(\kpp+2)(
u_y^2u_x^{-2\kpp-4}
-(2\kpp+3)^{-1}((2\kpp^2+5\kpp+4)u_yu_x^{-\kpp-2}-\kpp-1)
)
)\,U_3
\brk    
+(
(\kpp+2)^{-2}\,u_tu_x^{-2\kpp-3}((\kpp+1)^{-1}u_yu_x^{-\kpp-2}+1)
+u_y^2u_x^{-2\kpp-4}
\brk    
-2(2\kpp^2+7\kpp+1)(2\kpp+3)\,u_yu_x^{-\kpp-2}
+(2\kpp^2+8\kpp+7)(2\kpp+3)^{-1}
)\,U_4
\brk    
+u_tu_x^{-2\kpp-3}((\kpp+6)\,u_yu_x^{-\kpp-2}+(\kpp+1)(4\kpp^2+9\kpp+6))
\brk    
+(\kpp+1)^2(\kpp+2)(
(3\kpp+2)\,u_yu_x^{-2\kpp-4}
-(2\kpp+3)^{-1}(4(\kpp^2+3\kpp+3)u_yu_x^{-\kpp-2}+8\kpp^3
\brk    
+36\kpp^2+57\kpp+30)
)
)
\]
The structure equations are not involutive. The involutive system of structure equations includes equations for the differentials of the forms $\eta_1$, ... , $\eta_7$. These equations are too big to write them in full here.

We find contact integrable extensions of the form 
\[
d \omega = \left(
\sum \limits_{i=0}^3 A_i\,\theta_i
+\sum {}^{*} B_{ij}\,\theta_{ij}
+\sum \limits_{s=1}^7 C_s\,\eta_s
+\sum \limits_{j=1}^3 D_j\,\xi^j
+E\,\alpha
\right)\w\omega 
\]
\begin{equation}
\quad\quad\quad
+\sum \limits_{j=1}^3 \left(
\sum \limits_{k=0}^3 F_{jk}\,\theta_k + G_j\,\alpha
\right) \wedge \xi^j,
\label{simpliest_CIE}
\end{equation}
where  $\sum {}^{*}$ denotes suumation over all $i,j \in \mathbb{N}$ such that
$1\le i \le j \le 3$ and $(i,j)\not = (3,3)$. We consider two types of such extensions. The first one consists of extensions whose coefficients in right-hard side of  (\ref{simpliest_CIE}) depend on the invariants $U_1$, ... , $U_5$. 
The coefficients of extensions of the second type depend also on one additional function $W$ with the differential of the form
\begin{equation}
dW =\sum \limits_{i=0}^3 H_i \,\theta_i 
+ \sum {}^{*}   I_{ij}\,\theta_{ij}
+ \sum \limits_{s=1}^7 J_s\,\eta_s 
+ \sum \limits_{j=1}^3 K_j\,\xi^j 
+ \sum \limits_{q=0}^1 L_q\,\omega_q.
\label{dW_of_simpliest_CIE}
\end{equation}

We require Eqs. (\ref{SEs_of_rmmdKP}) and (\ref{simpliest_CIE}) or Eqs. (\ref{SEs_of_rmmdKP}), (\ref{simpliest_CIE}), and 
(\ref{dW_of_simpliest_CIE}) to be compatible. This con\-di\-ti\-on gives two contact integrable extensions of the first type defined by the formulas 
%
%
\[
\fl
d\omega_1=
((\kpp+2)^2\,(\alpha_1-(\kpp+1)\,\theta_0)-\theta_1-(\kpp+2)\,\theta_3)\w\xi^1
+\alpha_1\w\xi^2
\brk    
+((\kpp+2)\,(\alpha_1-(\kpp+1)\,\theta_0)-\theta_3)\w\xi^3
+(\alpha_1+\theta_2+\theta_{22}+
((\kpp+1)(\kpp\,U_1
\brk    
-(\kpp+2)\,U_3)-U_4
-(\kpp+1)(\kpp+2)(\kpp^2+6\kpp+4))(\kpp+1)^{-2}(\kpp+2)^{-2}\,\theta_0
\brk    
+((\kpp+1)((\kpp+1)\,U_2+(\kpp+2)\,U_3-\kpp(\kpp^2+3\kpp+3)\,U_1-\kpp(\kpp+1)(\kpp+2))
\brk    
+U_4)\,(\kpp+1)^{-2}\,\xi^1
+((\kpp+1)(\kpp^2\,U_1+(\kpp+2)\,U_3)-\kpp\,U_4
\]
\begin{equation}
\fl\qqq
-\kappa(\kappa+2)(3\kappa+2)(\kappa+1)^2)(\kpp+1)^{-2}(\kpp+2)^{-1}\,\xi^3
)\w\omega_1
\label{CIE_1_of_rmmdKP}
\end{equation}
and
%
%
\[
\fl
d\omega_2 = 
((\kappa+1)^2(\kappa+1)^2\,\theta_1-\theta_1+(\kappa+1)(\kappa+2)\,\theta_3)\w\xi^1
+\alpha_2\w\xi^2
\brk    
+((\kappa+1)(\kappa+1)\,\theta_2-\theta_3)\w\xi^3
+(\alpha_2+\theta_2+\theta_{22}+((\kpp+1)\,(\kpp\,U_1-(\kpp+2)\,U_3)
\brk    
-U_4
-(\kpp+1)(\kpp+2)(\kpp^2+6\kpp+4))(\kpp+1)^{-2}(\kpp+2)^{-2}\,\theta_0
\brk    
+(U_2-(\kpp+1)(\kpp+2)\,U_1)\,\xi^1
+(\kpp(\kpp+1)\,U_1 - U_4
\]
\begin{equation}
\fl\qqq
 -(\kpp+1)^2(\kpp+2)(3\kpp+2))(\kpp+1)^{-1}(\kpp+2)^{-2}\,\xi^3
)\w\omega_2
\label{CIE_2_of_rmmdKP}
\end{equation}
or  one contact integrable extension of the second type 
%
%
\[
\fl
d\omega_3 = 
((W+\kpp+2)^2\,\alpha_3-(W+\kpp+2)\,(\theta_{23}+(\kpp+1)(\kpp+2)\theta_0)-\theta_1)\w\xi^1
+\alpha_3\w\xi^2
\brk    
+((W+\kpp+2)\,\alpha_3-(\kpp+1)(\kpp+2)\,\theta_0-\theta_3)\w\xi^3
+(\alpha_3+\theta_2+\theta_{22}
+((\kpp+1)(\kpp\,U_1
\brk    
-(\kpp+2)\,U_3)-U_4-(\kpp+1)(\kpp+2)(\kpp^2+6\kpp+4))(\kpp+1)^{-2}(\kpp+2)^{-2}\,\theta_0
\brk    
-(((\kpp+1)(\kpp\,U_1-(\kpp+2)\,U_3)-U_4-(\kpp+1)(\kpp+2)(\kpp^2+6\kpp+4))\,W^2
\brk    
+(\kpp+2)((\kpp+1)\,((\kpp-1)\,U_1-2\,(\kpp+2)\,(U_3+(\kpp+1)(\kpp^2+5\kpp+3))-2\,U_4))\,W
\brk    
+(\kpp+2)^2((\kpp+1)(\kappa(\kappa^2+3\kappa+3)\,U_1-(\kappa+1)\,U_2-(\kappa+2)\,(U-\kappa^2(\kappa+1)))
\brk    
-U_4)(\kpp+1)^{-2}(\kpp+2)^{-2}\,\xi^1
-((\kpp+1)(\kappa\,U_1-(\kappa+2)\,(U_3+(\kappa+1)(\kappa^2+6\kappa+4)))
\]
\begin{equation}
\fl\qqq
-U_4)(\kpp+1)^{-2}(\kpp+2)^{-2}\,\xi^3
)\w\omega_3
\label{CIE_3_of_rmmdKP}
\end{equation}
\[
\fl
dW = -(\kappa+1)\,W\,(\alpha_3+\theta_0+\theta_2)
+Z\,\xi^2+(W+\kpp+2)(Z+(\kpp+1)\,W)\,\xi^3
\brk    
+(W+\kpp+2)((W+\kpp+2)\,Z+(\kpp+1)\,W\,(W-(\kpp+2)^{-1}\,U_1+3\kpp+4))\,\xi^1
\brk    
+(Z-(\kappa\,U_1-(\kappa+2)\,(U_3+(\kappa+1)^2(\kappa+6))
\]
\begin{equation}
\fl\qqq
-(\kpp+1)^{-1}U_4)(\kpp+2)^{-1}\,W)\,\omega_3
\label{dW_of_CIE_3_of_rmmdKP}
\end{equation}
with a parameter $Z$.

The inverse third  fundamental Lie theorem in Cartan's form,  \cite[\S 26]{Vasilieva1972}, \cite[p. 394]{Stormark2000}, 
guarantees existence of forms $\omega_1$, $\omega_2$, $\omega_3$  satisfying Eqs. (\ref{CIE_1_of_rmmdKP}), (\ref{CIE_2_of_rmmdKP}), and (\ref{CIE_3_of_rmmdKP}).  Since the forms $\theta_0$, ... , $\theta_{23}$, $\xi^1$, $\xi^2$, $\xi^3$ are known explicitly,  it is not hard to find the forms $\omega_i$. We have the following solutions to Eqs. (\ref{CIE_1_of_rmmdKP}), (\ref{CIE_2_of_rmmdKP}), and (\ref{CIE_3_of_rmmdKP}), respectively:   
%
%
\[
\fl
\omega_1 = 
\frac{u_{xx}}{u_xq_x}
\left(
dq 
-\left(\frac{u_t}{u_x} +(\kpp+2)\,\left(u_yu_x^\kpp+\frac{\kpp+1}{2\kpp+3}\,u_x^{2\kpp+2}\right)\right)\,q_x\,dt 
-q_x\,dx
\right.
\]
\begin{equation}
\fl
\qqq
\left.
-\left(\frac{u_y}{u_x}+u_x^{\kpp+1}\right)\,q_x\,dy
\right)
\label{WE_1_of_rmmdKP}
\end{equation}
%
%
\[
\fl
\omega_2 = 
\frac{u_{xx}}{u_xr_x}
\left(
dr 
-\left(\frac{u_t}{u_x} -(\kpp+1)(\kpp+2)\,\left(u_yu_x^\kpp-\frac{1}{2\kpp+3}\,u_x^{2\kpp+2}\right)\right)\,r_x\,dt 
-r_x\,dx
\right.
\]
\begin{equation}
\fl
\qqq
\left.
-\left(\frac{u_y}{u_x}-(\kpp+1)\,u_x^{\kpp+1}\right)\,r_x\,dy
\right)
\label{WE_2_of_rmmdKP}
\end{equation}
and
%
%
\[
\fl
\omega_3 = \frac{u_{xx}}{u_xs_x}\,\left(
ds -
\left(\frac{(\kappa+2)^2}{2\kappa+3}\,s_x^{2\kappa+3}-
(\kappa+2)\,\left(\frac{u_y}{u_x}+u_x^{\kappa+1}\right)\,s_x^{\kappa+2}
\right.
\right.
\brk    
\left.
\left.
+\left(
\frac{u_t}{u_x}+(\kappa+2)\,u_x^\kappa u_y+\frac{(\kappa+1)(\kappa+2)}{2\kappa+3}
u_x^{2\kappa+2}
\right)\,s_x\right)\,dt
-s_x\,dx
\right.
\]
\begin{equation}
\fl
\qqq
\left.
+\left(s_x^{\kappa+2} -\left(\frac{u_y}{u_x}+u_x^{\kappa+1}\right)\,s_x\right)\,dy
\right)
\label{WE_3_of_rmmdKP}
\end{equation}
with  $W=s_x^{\kappa+1}u_x^{-\kappa-1}$.

The forms (\ref{WE_1_of_rmmdKP}), (\ref{WE_2_of_rmmdKP}), (\ref{WE_3_of_rmmdKP})  are equal to zero if and only if the following over\-de\-ter\-mi\-ned systems of {\sc pde}s are satisfied: 
\begin{equation}
\fl
\left\{
\begin{array}{l}
q_t =\displaystyle{\left(\frac{u_t}{u_x} +(\kpp+2)\,\left(u_yu_x^\kpp+\frac{\kpp+1}{2\kpp+3}\,u_x^{2\kpp+2}\right)\right)\,q_x}
\\
q_y = \displaystyle{\left(\frac{u_y}{u_x}+u_x^{\kpp+1}\right)\,q_x}
\end{array}
\right.
\label{covering_1_of_rmmdKP}
\end{equation}
\begin{equation}
\fl
\left\{
\begin{array}{l}
r_t =\displaystyle{\left(\frac{u_t}{u_x} -(\kpp+1)(\kpp+2)\,\left(u_yu_x^\kpp-\frac{1}{2\kpp+3}\,u_x^{2\kpp+2}\right)\right)\,r_x}
\\
r_y = \displaystyle{\left(\frac{u_y}{u_x}-(\kpp+1)\,u_x^{\kpp+1}\right)\,r_x}
\end{array}
\right.
\label{covering_2_of_rmmdKP}
\end{equation}
\begin{equation}
\fl
\left\{
\begin{array}{l}
s_t =
\displaystyle{\frac{(\kappa+2)^2}{2\kappa+3}\,s_x^{2\kappa+3}-
(\kappa+2)\,\left(\frac{u_y}{u_x}+u_x^{\kappa+1}\right)\,s_x^{\kappa+2}}
\\
\qqq\qqq
+\displaystyle{\left(
\frac{u_t}{u_x}+(\kappa+2)\,u_x^\kappa u_y+\frac{(\kappa+1)(\kappa+2)}{2\kappa+3}
u_x^{2\kappa+2}
\right)\,s_x}
\\
s_y = \displaystyle{-s_x^{\kappa+2} +\left(\frac{u_y}{u_x}+u_x^{\kappa+1}\right)\,s_x}
\end{array}
\right.
\label{covering_3_of_rmmdKP}
\end{equation}
These systems are compatible whenever $u$ is a solution to Eq. (\ref{rmmdKP}), so these systems de\-fi\-ne differential coverings over (\ref{rmmdKP}). 

Expressing   $u_t$ and $u_y$ from (\ref{covering_1_of_rmmdKP}) and cross-differentiating yields 
\begin{equation}
q_{yy} = q_{tx}+\left((\kpp+1)\,\frac{q_y^2}{q_x^2}-\frac{q_t}{q_x}\right)\,q_{xx}-\kpp\,\frac{q_y}{q_x}\,q_{xy}
\label{rmmdKP-I}
\end{equation}
Previously  Eq. (\ref{rmmdKP-I}) and the B\"acklund transformation (\ref{covering_1_of_rmmdKP}) were found in \cite{Morozov2010} by means of another method.

From Eqs. (\ref{covering_2_of_rmmdKP}) we have
\begin{equation}
\fl
\left\{
\begin{array}{l}
u_t =\displaystyle{\left(\frac{r_t}{r_x} 
+(\kpp+1)(\kpp+2)\,\left(\frac{r_y}{r_x}\,u_x^{\kpp+1}+\frac{(\kpp+2)(2\kpp+1)}{2\kpp+3}\,u_x^{2\kpp+2}\right)\right)\,u_x}
\\
u_y = \displaystyle{\left(\frac{r_y}{r_x}+(\kpp+1)\,u_x^{\kpp+1}\right)\,u_x}
\end{array}
\right.
\label{inverse_BT_2_for_rmmdKP}
\end{equation}
The compatibility condition for this system is
\[
\fl 
(u_t)_y-(u_y)_t = -(\kpp+1)(\kpp+2)\,u_x^{\kpp+2}r_x^{-2}\,\left(G\,r_x - \kpp\,(\kpp+2)\,u_x^{\kpp+1} \,(r_y\,r_{xx}-r_x\,r_{xy})\right)
=
\]
\begin{equation}
\fl\qqq\qqq
=0
\label{compatibility_of_inverse_BT_2_for_rmmdKP}
\end{equation}
where
\[
G = 
r_{yy} - r_{tx}-\left((\kpp+1)\,\frac{r_y^2}{r_x^2}-\frac{r_t}{r_x}\right)\,r_{xx}+\kpp\,\frac{r_y}{r_x}\,r_{xy}
\]
When  $\kappa=0$, system (\ref{inverse_BT_2_for_rmmdKP}) is compatible whenever $G =0$, that is, whenever $r$ is a solution to Eq. (\ref{rmmdKP-I}).  When $\kappa\not= 0$, Eq. (\ref{compatibility_of_inverse_BT_2_for_rmmdKP}) entails $u_x^{\kpp+1} = H$ with
\[
H = - \kpp^{-1}(\kpp+2)^{-2}\,G\,r_x\,(r_y\,r_{xx}-r_x\,r_{xy})^{-1}
\]
Substituting this into (\ref{inverse_BT_2_for_rmmdKP}) gives a system of {\sc pde}s with the compatibility condition
\[
\fl
\kpp\,(2\kpp+3)\,r_x^2\,H_t
-\kpp\,(\kpp+2)\,r_x\,\left(2(\kpp+2)(2\kpp+1)\,r_xH+(2\kpp+3)\,r_y\right)\,H_y
\brk    
+\kpp\,\left(
(\kpp+1)(\kpp+2)^2(2\kpp+1)\,r_x^2\,H^2+2(\kpp+2)(2\kpp+1)\,r_xr_yH
\right.
\brk    
\left.
-(2\kpp+3)(r_tr_x-(\kpp+2)\,r_y^2)
\right)\,H_x
-(\kpp+1)\,\left(
(2\kpp^2+5\kpp+1)\,r_xG
\right.
\]
\begin{equation}
\left.
+\kpp(2\kpp+3)(r_xr_{tx}-r_tr_{xx})
\right)\,H
-(2\kpp+3)\,r_y\,G =0
\label{covering_equation_2_for_rmmdKP}
\end{equation}
Thus Eqs. (\ref{covering_2_of_rmmdKP}) define a B\"acklund transformation from Eq. (\ref{rmmdKP}) to 
the third order equa\-ti\-on (\ref{covering_equation_2_for_rmmdKP}) for $r$.

Finally, excluding $u$ from (\ref{covering_3_of_rmmdKP}) shows that $s$ is a solution to the same equation (\ref{rmmdKP}). So,
(\ref{covering_3_of_rmmdKP}) defines an auto-B\"acklund transformation for Eq. (\ref{rmmdKP}). This transformation was found in 
\cite{MorozovPavlov2010}

\section*{Acknowledgements}

I am very grateful to M.V. Pavlov for many stimulating discussions.

\end{document}